\documentclass[12pt,leqno]{article}
\usepackage{amsfonts}
\usepackage{amssymb}
\usepackage{enumerate}
\usepackage{amsmath}
\usepackage[T1]{fontenc}

\newtheorem{ttt}{Theorem}[section]
\newtheorem{llll}[ttt]{Lemma}
\newtheorem{ccc}[ttt]{Claim}
\newtheorem{eee}[ttt]{Example}
\newtheorem{fff}[ttt]{Fact}
\newtheorem{rrr}[ttt]{Remark}
\newtheorem{sss}[ttt]{Statement}
\newtheorem{ddd}[ttt]{Definition}
\newtheorem{qqq}[ttt]{Question}
\newtheorem{cccc}[ttt]{Corollary}
\newtheorem{nnn}[ttt]{Notation}

\newcommand{\bt}{\begin{ttt}}
\newcommand{\bl}{\begin{llll}}
\newcommand{\bc}{\begin{ccc}}
\newcommand{\bex}{\begin{eee}}
\newcommand{\bfa}{\begin{fff}}
\newcommand{\br}{\begin{rrr}\upshape}
\newcommand{\bs}{\begin{sss}}
\newcommand{\bd}{\begin{ddd}\upshape}
\newcommand{\bq}{\begin{qqq}}
\newcommand{\bno}{\begin{nnn}}
\newcommand{\bcor}{\begin{cccc}}

\newcommand{\bp}{\noindent\textbf{Proof. }}

\newcommand{\et}{\end{ttt}}
\newcommand{\el}{\end{llll}}
\newcommand{\ec}{\end{ccc}}
\newcommand{\eex}{\end{eee}}
\newcommand{\efa}{\end{fff}}
\newcommand{\er}{\end{rrr}}
\newcommand{\es}{\end{sss}}
\newcommand{\ed}{\end{ddd}}
\newcommand{\eq}{\end{qqq}}
\newcommand{\ecor}{\end{cccc}}
\newcommand{\eno}{\end{nnn}}

\newcommand{\ep}{\hspace{\stretch{1}}$\square$\medskip}

\newcommand{\lab}[1]{\label{#1}}

\newcommand{\NN}{\mathbb{N}}
\newcommand{\ZZ}{\mathbb{Z}}
\newcommand{\QQ}{\mathbb{Q}}
\newcommand{\RR}{\mathbb{R}}
\newcommand{\CC}{\mathbb{C}}
\newcommand{\TT}{\mathbb{T}}

\newcommand{\al}{\alpha}

\newcommand{\om}{\omega}
\newcommand{\si}{\sigma}
\newcommand{\ka}{\kappa}

\newcommand{\iH}{\mathcal{H}}

\newcommand{\iS}{\mathcal{S}}
\newcommand{\iN}{\mathcal{N}}

\def\su{\subset}
\def\sm{\setminus}

\newcommand{\cof}{{\rm cof}}
\newcommand{\cov}{{\rm cov}}

\title{Covering locally compact groups by less than $2^\om$ many translates 
of a compact nullset} 

\author{M\'arton Elekes\thanks{Partially supported by Hungarian Scientific
Foundation grants no.~49786, 37758 and F 43620.} \ and \'Arp\'ad
T\'oth\thanks{Partially supported by Hungarian Scientific Foundation
grant no.~T 049693.}}

\begin{document}

\maketitle 

\begin{abstract}
Gruenhage asked if it was possible to cover the real line by less than 
continuum many translates of a compact
nullset. Under the Continuum Hypothesis the answer is obviously negative. 
Elekes and Stepr\=ans \cite{ES} gave an affirmative answer 
by showing that if $C_{EK}$ is the well known
compact nullset considered first by Erd\H os and Kakutani \cite{EK} then 
$\RR$ can be covered by $\cof(\iN)$ many translates of $C_{EK}$.
As this set has no analogue in more general groups, it was posed as an
open question in \cite{ES} whether such a result holds for uncountable 
locally compact Polish groups. In this paper we give an affirmative answer in
the abelian case. 

More precisely, we show that if $G$ is a nondiscrete locally compact
abelian group in which every open subgroup is of index at most $\cof(\iN)$
then there exists a compact set $C$ of Haar measure zero such that $G$ 
can be covered by $\cof(\iN)$ many translates of $C$. This result, which is
optimal in a sense, covers the
cases of uncountable compact abelian groups and of nondiscrete separable
locally compact abelian groups.

We use Pontryagin's duality theory to reduce the problem to three
special cases; the circle group, countable products of finite discrete abelian
groups, and the groups of $p$-adic integers, and then we solve the problem on
these three groups separately.

In addition, using representation theory, we reduce the nonabelian case to the
classes of Lie groups 
and profinite groups, and we also settle the problem for Lie 
groups.\footnote{M.~Ab\'ert recently gave an affirmative answer for profinite
groups \cite{Ab}.} 
\end{abstract}

\insert\footins{\footnotesize{MSC codes: Primary 22B05, 28C10, 03E17; Secondary
28E15, 03E35}}
\insert\footins{\footnotesize{Key words: compact, null, zero, translate, 
continuum, $2^\om$, Pontryagin duality, LCA, Polish, separable, group, 
cof, consistent, ZFC, $p$-adic integers, profinite groups, Lie groups}}

\section{Introduction}

Under the Continuum Hypothesis the real line
obviously cannot be covered by less than $2^\omega$ many translates of a set
of Lebesgue measure zero. On the other hand, it is well known that in some 
models of set theory there
exists such a covering \cite{BJ}. Moreover, we can obviously assume that
the set is $G_\delta$. Gruenhage \cite{Gr} asked whether such a covering can be
constructed with an $F_\sigma$ or closed or compact nullset (using of course
some extra set-theoretic assumption).

\bq \lab{q:gruenhage}
(Gruenhage) Let $C\subset
\RR$ be a compact set of Lebesgue measure zero and $A\subset \RR$ be of
cardinality less than $2^\omega$. Does that imply $C+A\not=\RR$?
\eq

We remark that it is well known that 
in some models of set theory $\RR$ can be covered by less than
$2^\om$ many compact nullsets (\cite{BJ} or \cite{BS}), but in these
coverings the sets are not translates of each other.

We also remark that already \cite{Mi} considers cardinal invariants of closed
measure zero sets, and \cite{MS}, \cite{Pa} and \cite{Sh} deals with 'translative
cardinal invariants'; that is, when the small sets considered are translates
of each other.
For another very closely related paper see \cite{Zi}.

Gruenhage gave an affirmative answer to Question \ref{q:gruenhage} when $C$ is the
classical Cantor set \cite{Gr}, and later Darji and Keleti \cite{DK}
generalized his results to the class of compact nullsets of so called packing
dimension less than 1.

Then Elekes and Stepr\=ans \cite{ES} answered all versions of
Gruenhage's question in the negative as follows. 

\bd 
Denote
\[
C_{EK} = \left\{ \left. \sum_{n=2}^{\infty} \frac{d_n}{n!} \ \right| \ d_n
\in\{0,1,\dots, n-2\}\ \forall n \right\}.
\]
\ed

The letters E and K stand for Erd\H os and Kakutani.

\bd
Let $\iN$ denote the set of Lebesgue nullsets of the real line, and let
$\cof(\iN) = \min\{|\iH| : \iH \su \iN, \forall N\in\iN \ \exists H\in\iH, N\su
H\}$. 
\ed

It is not hard to see that $\om < \cof(\iN) \le 2^\om$ \cite{BJ}.

\bt\lab{t:ES}
(Elekes-Stepr\=ans) $\RR$ can be covered by $\cof(\iN)$ many
translates of the compact nullset $C_{EK}$. 
\et

As $\cof(\iN) < 2^\om$ is consistent with the axioms of set theory \cite{BJ},
we obtain the following.

\bcor\lab{covers}
It is consistent with the axioms of set theory that less than continuum many
translates of a compact set of measure zero cover the real line.
\ecor

As $C_{EK}$ has no analogue in more general groups, it was posed as an
open question in \cite{ES} whether such a result holds for uncountable 
locally compact Polish groups. The main goal of this paper is to show 
that the answer is affirmative in the abelian case 
(Corollary \ref{c:main}). Note that countable locally compact groups
are not interesting from this viewpoint, and also that the assumption that the
group is Polish is also natural, since our problem actually considers a cardinal
invariant (see \cite{BJ}), and this topic is usually discussed  in the
framework of Polish spaces. 

First we use Pontryagin's duality theory to reduce the problem to three special
cases; the circle group, countable products of finite discrete groups,
and the groups of $p$-adic integers, then we solve the problem separately for
these groups.

In Section \ref{s:nonabelian} we discuss the nonabelian version of our
problem. We reduce the nonabelian case to the cases of
Lie groups and profinite groups, and we show that every nondiscrete 
Lie group in which every open subgroup is of index at most $\cof(\iN)$
can be covered by $\cof(\iN)$ many
left translates of a compact set of Haar measure zero.

Note that a set is of left Haar measure zero iff it is of right Haar measure
zero.

All groups are tacitly assumed to be Hausdorff.

\br
The referee pointed out the following interesting facts.

\noindent
1. Our method of reducing the problem to some special groups is
   fairly general. Therefore it may well be applicable to show that all locally
   compact groups possess a certain property, supposing that whenever a factor
   group $G/H$ has the property then $G$ itself has it as well. 

\noindent
2. The use of $\cof(\iN)$ is not optimal, one can show that
   consistently it can be improved. In fact, it could be replaced with the
   the least cardinality $\kappa$ for which for every pair $f,g:\om\to\om$
   converging to infinity every $f$-slalom can be covered by $\kappa$ many
   $g$-slaloms (see Definition \ref{d:slalom}). However, as this is not a very
   well known invariant (there is 
   no abbreviation for it), and most probably this is also not optimal, we
   still prefer to use $\cof(\iN)$. 
 
\noindent
3. Question \ref{q:gruenhage} is closely related to the following, which
   essentially asks whether the set of translations we use can be arbitrary. 
   \emph{Is it true
   that for every uncountable $X\su\RR$ there exists a countable set $Y$ and a
   closed nullset $F$ such that $(X+Y)+F=\RR$?} On can easily show that this is
   in fact equivalent to the following. \emph{Is it true that for every
   uncountable 
   $X\su\RR$ there exists an $F_\si$ nullset $A$ so that $X+A=\RR$?} On can
   very easily give a consistent negative answer to these questions (e.g. if
   $\cov(\iN)=2^\om > \om_1$), but a negative answer in $ZFC$ would be
   interesting. On the other hand, a consistent affirmative answer would prove
   the consistency of the so called Borel Conjecture + Dual
   Borel Conjecture, which is a longstanding open problem.
\er

\section{The abelian case}

\br
It may be instructive to bear in mind that the proof (just as in Section
\ref{s:nonabelian} in the nonabelian case) will consist of two parts. First
we prove a purely \emph{analysis} result constructing a compact nullset and
showing that every so called `slalom' can be covered by a translate of the
set, then we apply a purely \emph{set-theoretic} result stating that
consistently less than $2^\om$ many slaloms can cover the space.
\er

A topological group is LCA if it is locally compact and abelian. 

\bd
We say that a locally compact group $G$ is \emph{nice} if 
there exists a compact set $C \su G$ of Haar
measure zero such that $G$ can be covered by $\cof(\iN)$ many
left translates of $C$.
\ed

The aim of this section is to prove the following.

\bt\lab{main}
Suppose that $G$ is a nondiscrete LCA group in which every open 
subgroup is of index at most $\cof(\iN)$. Then $G$ is nice; that is, there
exists a compact set $C\su G$ of Haar measure zero such that $G$ can be
covered by $\cof(\iN)$ many translates of $C$.
\et

\br
Both conditions of the theorem are necessary. First, if $G$ is
discrete then the only nullset is the emptyset, so no covering by
nullsets exists. Secondly, if there is an open 
subgroup of index $\kappa$ then at least $\kappa$ many compact nullsets
are needed to cover $G$, since a compact set can only intersect
finitely many cosets. 

In fact, as ``$\cof(\iN)=\omega_1$ and $2^\omega =\omega_2$'' is consistent
with the axioms of set theory \cite{BJ}, we actually obtain
the following consistent characterisation.

\bcor\lab{c:char}
It is consistent with the axioms of set theory that an LCA group $G$ can be
covered by 
less than $2^\om$ many translates of a compact nullset iff $G$ is
nondiscrete and has no open subgroup of index at least $2^\om$. 
\ecor

\er

Before the proof of Theorem \ref{main} we formulate two more corollaries.

\bcor
Every uncountable compact abelian group and every nondiscrete separable LCA
group is nice; that is, it can be covered by $\cof(\iN)$ many translates of a
compact nullset.
\ecor

As $\cof(\iN) < 2^\om$ is consistent with the axioms of set theory \cite{BJ},
and as every Polish space is separable, we obtain the following, which answers
Question 3.2 in \cite{ES} in the abelian case.

\bcor\lab{c:main}
It is consistent with the axioms of set theory that every uncountable locally
compact abelian Polish group can be covered by less than $2^\om$ 
many translates of a compact nullset. 
\ecor

In the rest of the section we prove Theorem \ref{main}. First we need two
technical lemmas.



\bl\lab{l:coord}
Let $n\ge 0$ be an integer, $G$ be a finite group, $A$ and $S$ be
subsets of $G$ 
such that $\left(1-\frac{1}{n+3}\right)|G| \le |A|$ and $|S| \le n+2$. 
Then there exists $g\in G$ such that $S \su gA$. 
\el

\bp
Clearly $S \not\su gA$ iff $g \in S(G\sm A)^{-1}$. So it is enough to
check that $S(G\sm A)^{-1} \neq G$, which is clear, since 
\[
|S(G\sm A)^{-1}| \le |S| \cdot |(G\sm A)| \le (n+2)\frac{|G|}{n+3}  < |G|.
\]
\ep

For a sequence $(X_n)_{n\in\NN}$ of set $\times_{n\in\NN} X_n$ denotes the
Cartesian product. 

\bd\lab{d:slalom}
For every $n\in\NN$ let $X_n$ be an arbitrary set, and fix a function
$f:\NN\rightarrow\NN\setminus\{0\}$. An \emph{$f$-slalom} is a set of the form
\[
S=\times_{n\in\NN} S_n, \textrm{ where } S_n\subset X_n, \ |S_n|\le f(n) \
\ (n\in\NN).
\]
\ed

\bl\lab{l:slalom}
Let $f_0:\NN\to\NN\sm\{0\}$ be such that $\lim_\infty f_0 = \infty$, 
and let $X_n$ $(n\in\NN)$ be countable
sets. Then $\times_{n\in\NN} X_n$ can be covered by $\cof(\iN)$ many
$f_0$-slaloms. 
\el

\bp
\cite[2.3.9]{BJ} states that there exist a system of functions $f_\al
: \NN \to \NN\sm\{0\}$ $(\al<\cof(\iN))$ with $\sum_{n\in\NN^+}
{f_\al(n)}/{n^2} < \infty$, and for every $\al<\cof(\iN)$
there exists an $f_\al$-slalom $S_\al = \times_{n\in\NN} (S_\al)_n \su
\NN^\NN$ such that these slaloms cover $\NN^\NN$
 mod finite; that is, for every $g\in\NN^\NN$ there exists
$\al<\cof(\iN)$ such that $\{n\in\NN : g(n) \notin (S_\al)_n\}$ is finite.
For an $f$-slalom $S\su \NN^\NN$ let $\iS_S = \{S' \su \NN^\NN : S'
\textrm{ is an $f$-slalom, and } \{n\in\NN : S_n \neq S_n'\}
\textrm{ is finite}\}$. Clearly, every $\iS_S$ is countable, and hence
$\cup_{\al<\cof(\iN)} \iS_{S_\al}$ is easily seen to be a set of
$\cof(\iN)$ many slaloms actually covering $\NN^\NN$. So we can assume
that $\cup_{\al<\cof(\iN)} S_\al = \NN^\NN$. Put $f(n) = n^2+1$. 
Clearly, $\{n\in\NN : f_\al(n) > f(n)\}$ is finite for every $\al$,
and therefore an argument similar to the previous one yields that
every $f_\al$-slalom can be covered by countably many $f$-slaloms. So
$\NN^\NN$ can be covered by $\cof(\iN)$ many $f$-slaloms.

\cite[2.10.]{GL} states that if $f,g : \NN \to \NN\sm\{0\}$ are such
that $\lim_\infty f = \lim_\infty g = \infty$, then the minimal number of
$f$-slaloms needed to cover $\NN^\NN$ equals the minimal number of
$g$-slaloms needed to cover $\NN^\NN$. Therefore
$\NN^\NN$ can be covered by $\cof(\iN)$ many $f_0$-slaloms, hence
$\times_{n\in\NN} X_n$ can also be covered by $\cof(\iN)$ many
$f_0$-slaloms.
\ep

In order to prove Theorem \ref{main} first we need to prove it in two
special cases; for countable products of finite discrete (abelian)
groups and for the groups of $p$-adic integers.

For a sequence $(G_n)_{n\in\NN}$ of compact groups
$\otimes_{n\in\NN} G_n$ is the (Cartesian) product group with the product
topology.

\bt\lab{t:product}
For every $n\in\NN$ let $G_n$ be a discrete finite group of at least
$2$ elements. Then $\otimes_{n\in\NN} G_n$ is nice.
\et

\bp
Write $\NN$ as the disjoint union of finite sets $N_n$ such that
$2^{|N_n|} > 2(n+3)$, and define $G_n' = \otimes_{k\in N_n} G_k$. Then 
$\otimes_{n\in\NN} G_n = \otimes_{n\in\NN} G_n'$ and $|G_n'| >
2(n+3)$. Hence for every $n\in\NN$ we can find an $A_n\su G_n'$ such that
\[
\left(1-\frac{1}{n+3}\right)|G_n'| \le |A_n| \le
\left(1-\frac{1}{2(n+3)}\right)|G_n'|.
\]
Define $C = \times_{n\in\NN}
A_n$. Then $C$ is clearly compact, and $\Pi_{n\in\NN}
\left(1-\frac{1}{2(n+3)}\right) = 0$ implies that $C$ is of Haar measure
zero. 

Put $f_0(n) = n+2$ $(n\in\NN)$. By Lemma \ref{l:slalom} $\otimes_{n\in\NN}
G_n'$ can be covered by $\cof(\iN)$ many $f_0$-slaloms. 
We will complete the proof by showing that every $f_0$-slalom $S =
\times_{n\in\NN} S_n \su \otimes_{n\in\NN} G_n'$ can be covered by a
left translate of $C$. 

For every $n\in\NN$ we can apply Lemma \ref{l:coord}
to $n, G_n', A_n$ and $S_n$, and so we obtain a $g_n \in G_n'$ such
that $S_n \su g_n A_n$. But then for $g = (g_n)_{n\in\NN} \in
\otimes_{n\in\NN} G_n'$ we have $S \su \times_{n\in\NN} \ g A_n =
g C$.
\ep

We need certain properties of the $p$-adic integers $\ZZ_p$ that we collect
here for the convenience of the reader. For a precise treatment see
e.g.~\cite {Ro}. 
The underlying topological space is
$\{0,1,\dots,p-1\}^\NN$ equipped with the product topology (each factor is considered
discrete). Addition is coordinatewise with carried digits from the
$n^{th}$ coordinate to the $n+1^{st}$; that is, if $x=(x_n)_{n\in\NN},
y=(y_n)_{n\in\NN} \in \ZZ_p$ then 
$(x+y)_0 = x_0 + y_0$ if $x_0 + y_0 \le p-1$ while $(x+y)_0 =
x_0 + y_0 - p$ if $x_0 + y_0 \ge p$. In the second case when
calculating $(x+y)_1$ we add $1$ to $x_1+y_1$ and then check whether
the sum is greater than $p-1$, etc, recursively.

%
%
%
%

\bt\lab{t:p-adic}
For every prime $p$ the group of $p$-adic integers $\ZZ_p$ is nice.
\et

\bp
If we forget about the group operation then we can write
$\ZZ_p = \times_{n\in\NN} X_n$, where $X_n = \{0,1,\dots,p-1\}$ 
for every $n\in\NN$.

Write $\NN$ as the disjoint union of the finite \emph{intervals}
$[k_n, k_{n+1})$, where $\{k_n\}_{n\in\NN}$ is a strictly monotone
increasing sequence of nonnegative integers such that  
$p^{k_{n+1}-k_n} > 2(n+3)$.
Define $X_n' = \times_{k\in [k_n, k_{n+1})} 
X_k$. Then 
$\times_{n\in\NN} X_n = \times_{n\in\NN} X_n'$ and $|X_n'| >
2(n+3)$. As above, for every $n\in\NN$ we can find an $A_n\su X_n'$ such that
\[
\left(1-\frac{1}{n+3}\right)|X_n'| \le |A_n| \le
\left(1-\frac{1}{2(n+3)}\right)|X_n'|.
\]
Let $C = \times_{n\in\NN}
A_n$. Again, $C$ is compact and is of Haar measure zero. 

Put $f_0(n) = \lfloor \frac{n+2}{2} \rfloor$ $(n\in\NN)$.
($\lfloor x \rfloor$ is the integer part of $x$.) 
By Lemma \ref{l:slalom} $\times_{n\in\NN}
X_n'$ can be covered by $\cof(\iN)$ many $f_0$-slaloms. 
We will complete the proof by showing that every $f_0$-slalom $S =
\times_{n\in\NN} S_n \su \times_{n\in\NN} X_n'$ can be covered by a
translate of $C$.

For every $n\in\NN$ we define a new group $G_n$ 
(\emph{not} a subgroup of $\ZZ_p$) as follows. 
Let $G_n = X_n' = \times_{k\in [k_n,k_{n+1})} X_k$, and for $x = (x_k)_{k\in [k_n,k_{n+1})} \in
G_n$ and $y = (y_k)_{k\in [k_n,k_{n+1})} \in G_n$ put
\[
(x +_{G_n} y)_k =  (x +_{\ZZ_p} y)_k \textrm{ for every } k \in [k_n,k_{n+1});
\]
that is, we always forget about the last carried digit. One can check that
$G_n$ with this addition is indeed a group. For example, to avoid all
calculations, it is easy to see that this group is
(canonically isomorphic to) $p^{k_n} \ZZ_p / p^{k_{n+1}} \ZZ_p$ and
also to $p^{k_n} \ZZ / p^{k_{n+1}} \ZZ$, but we will not use this fact.

Put $1_n = \chi_{\{k_n\}}$. ($\chi_H$ is the characteristic function of the
set $H$.) Fix $n\in\NN$, and set $\tilde{S}_n = S_n \cup (S_n +_{G_n} 1_n)$.
As $|S_n| \le \lfloor \frac{n+2}{2} \rfloor$, we
clearly have $|\tilde{S}_n| \le n+2$, hence we can apply Lemma
\ref{l:coord} to $n$, $G_n$, $A_n$ and $\tilde{S}_n$, and so we
obtain a $g_n \in G_n = \times_{k\in [k_n,k_{n+1})} X_k$ such that $\tilde{S}_n \su A_n
+_{G_n} g_n$. Let $x_n$ be the inverse of $g_n$ in $G_n$, then $\tilde{S}_n
+_{G_n} x_n \su A_n$. Put $x =
(x_n)_{n\in\NN} \in \times_{n\in\NN} X_n'$. We claim that $S +_{\ZZ_p} x \su
C$, which will complete the proof. Fix $s = (s_n)_{n\in\NN} \in S$. When we
recursively calculate the digits of $s +_{\ZZ_p} x$, we need to show that for
every $n\in\NN$ we have $\left( (s +_{\ZZ_p} x)_k \right)_{k\in [k_n,k_{n+1})} \in A_n$,
but this is clear, as $\left( (s +_{\ZZ_p} x)_k \right)_{k\in [k_n,k_{n+1})}$ equals
either $\left( \left(s_n +_{G_n} x_n \right)_k
\right)_{k\in [k_n,k_{n+1})}$ or $\left( s_n +_{G_N} x_n +_{G_n}
1_n)_k \right)_{k\in [k_n,k_{n+1})}$, depending on whether there is a
carried digit at $k_n$ or not. 
\ep

Before proving Theorem \ref{main} we need an 
algebraic fact about abelian groups. It is formulated in Theorem \ref{t:group},
which is well known, e.g.~a more general version
appears in \cite{KR}, but for the sake of completeness we include a proof
below.

\bd
Let $G$ be a group. 
For every
$n\in\NN$ let $G_{p^n} = \{g\in G : p^n g =0 \}$, and also let $G_{p^\infty} =
\cup_{n\in\NN} G_{p^n}$. We say that $G$ is a \emph{$p$-group} if 
$G = G_{p^\infty}$.
\ed

\bd
Let $p$ be a prime. An abelian group $G$ is called \emph{quasicyclic} if it is
generated by a sequence $(g_n)_{n\in\NN}$ with the property that $g_0 \neq
0$ and $p g_{n+1} = g_n$ for every $n\in\NN$. For a fixed prime $p$ the unique
(up to isomorphism) quasicyclic group is denoted by $C_{p^\infty}$.
\ed

Note that $C_{p^\infty} = (\QQ/\ZZ)_{p^\infty} = (\RR/\ZZ)_{p^\infty}$.


\bl\lab{l:quasi}
Let $p$ be a prime and $G$ be an infinite abelian $p$-group such that
$G_{p^n}$ is finite for every $n\in\NN$. Then $G$ contains
$C_{p^\infty}$ as a subgroup.
\el

\bp
We define a graph on $G$ as follows. For every nonzero $g\in G$ we connect
$pg$ with $g$. The obtained graph is clearly a tree (with root $0$) in which
each node has finitely many immediate successors by the finiteness of the
$G_{p^n}$'s. So by K\"onig's lemma \cite[5.7]{Ku} the tree has an infinite
branch, which clearly generates a quasicyclic subgroup.
\ep

For a sequence $(G_n)_{n\in\NN}$ of abelian groups $\oplus_{n\in\NN} G_n$ is
the direct sum group (that is, those elements of the product that only have
finitely many nonzero coordinates) with the discrete topology.

\bt\lab{t:group}
Every infinite abelian group $G$ contains a subgroup isomorphic to one of the
following. 
\begin{enumerate}
\item $\ZZ$,
\item $\oplus_{n\in\NN} G_n$, where each $G_n$ is a finite abelian group of 
at least two elements,
\item $C_{p^\infty}$ for some prime $p$.
\end{enumerate}
\et

\bp
If $G$ contains an element of infinite order then $G$ contains
$\ZZ$ as a subgroup. Therefore we may assume that $G$ is a torsion group.

Every torsion group is the direct sum of $p$-groups; $G =
\oplus_{p \textrm { prime }} G_{p^\infty}$ \cite[2.1]{Fu}. 

Suppose that $|G_{p^\infty}| \ge 2$ for infinitely many primes
$p$. For every such $p$ we can find a finite nontrivial subgroup of
$G_{p^\infty}$, and hence we have a sequence $(G_n)_{n\in\NN}$ of finite
nontrivial groups such that $\oplus_{n\in\NN} G_n \su G$. 
So we may assume that $|G_{p^\infty}| = 1$ for all but finitely many
primes. As $G$ is infinite, there is a prime $p$ for which $G_{p^\infty}$ is
infinite.

Assume that $G_p$ is infinite. Then $G_p$ is clearly an infinite
dimensional vector field over $\mathbb{F}_p$, therefore it contains 
$\oplus_{n\in\NN} C_p$ as a subgroup. ($C_p$ is the cyclic group of $p$
elements.)

So we may assume that $G_p$ is finite. Then we claim that
$G_{p^n}$ is also finite for every $n\in\NN$. We prove this by induction
on $n$. The map $g \mapsto pg$ is a homomorphism of $G_{p^{n+1}}$ into
$G_{p^n}$ with kernel $G_p$, so $|G_{p^{n+1}}| \le
|G_{p^n}| \cdot |G_p|$, which finishes the induction. Hence we can
apply Lemma \ref{l:quasi} to $G_{p^\infty}$, and we obtain that
$C_{p^\infty} \su G_{p^\infty} \su G$. This finishes the proof.
\ep

The following lemma is crucial.

\bl\lab{l:factor}
Let $G$ be a locally compact group and $H\su G$ a compact normal subgroup. If
$G/H$ is nice then $G$ is also nice.
\el

\bp
Let $\mu_G$ be a left Haar measure on $G$, and let
$\pi:G\to G/H$ be the canonical homomorphism. Then by 
\cite[\S 63.~Thm.~C.]{Ha} $\mu_G \circ\pi^{-1}$ is a left Haar measure
on $G/H$. This shows that the inverse image of a nullset in $G/H$ under $\pi$
is a nullset in $G$. Moreover, \cite[\S 63.~Thm.~B.]{Ha} states that
the inverse image of a compact set under $\pi$ is also compact.

Hence if $C\su G/H$ is a compact nullset witnessing that $G/H$ is
nice then $\pi^{-1}(C) \su G$ is a compact nullset witnessing that 
$G$ is also nice. 
\ep

\br
The following example shows that the lemma does not hold in general,
that is, when $H$ is a closed normal subgroup. Let $H$ be a discrete
group of cardinality greater than $\cof(\iN)$ and let $G = H \times
\RR$. Then $G/H$ is nice by Theorem
\ref{t:ES}, but $G$ is not nice as every compact set intersects 
only finitely many cosets.
\er

Now we are ready to prove our main theorem.

\bigskip

\bp (Theorem \ref{main}) By the Principal Structure Theorem of LCA Groups
\cite[2.4.1]{Ru}, $G$
has an open subgroup $H$ which is of the form $H = K \otimes \RR^n$, where $K$
is a compact subgroup and $n\in\NN$. By assumption the
index of $H$ is at most $\cof(\iN)$, so it suffices to prove that 
$H$ is nice, so we can assume $G=H$. 

Suppose $n\ge 1$. By \cite[2.1]{ES} $\RR$ is nice, let $C$ be the compact
nullset witnessing this fact. Then it is easy to see that $K \times C \times
[0,1]^{n-1}$ witnesses that $G = K \otimes \RR^n$ is nice. Hence we can assume
$n=0$, so $G$ is compact.

By Lemma \ref{l:factor} it is sufficient to find a closed subgroup $H\su G$
such that $G/H$ is nice. By \cite[2.1.2]{Ru} (and by the Pontryagin
Duality Theorem \cite[1.7.2]{Ru}) factors of $G$ are the same as
(isomorphically homeomorphic to) dual groups of closed subgroups of $\hat{G}$. 
As $G$ is compact, $\hat{G}$ is discrete \cite[1.2.5]{Ru}. Hence it suffices
to find a subgroup $M\su \hat{G}$ such that $\hat{M}$ is nice. 

By Theorem \ref{t:group} $\hat{G}$ has a subgroup isomorphic either to $\ZZ$,
or to $\oplus_{n\in\NN} G_n$, where each $G_n$ is a finite abelian group of 
at least two elements, or to $C_{p^\infty}$ for some prime $p$. We need to
show that the duals of these groups are nice.

By \cite[2.1]{ES} $\RR$ is nice, which easily implies
that the circle group $\TT$ is also nice, so we are done in the first case, 
since $\hat{\ZZ} = \TT$. 

In the second case note that $\hat{G}$ is finite iff $G$ is finite, hence each
$\hat{G}_n$ is finite. 
By \cite[2.2.3]{Ru} the dual of a direct sum (equipped with the discrete
topology) is the direct product of the dual groups (equipped with the product
topology), so $(\oplus_{n\in\NN} G_n)\hat{} = \otimes_{n\in\NN} \hat{G}_n$,
which is nice by Theorem \ref{t:product}.  

Finally, the third case is settled by Theorem \ref{t:p-adic}, since by
\cite[25.2]{HR} $\hat{C}_{p^\infty} = \ZZ_p$.
\ep




\section{The nonabelian case}\lab{s:nonabelian}

The aim of this section is to reduce the general case to the case of
profinite groups; that is, inverse limits of finite
groups.\footnote{We have been informed by M.~Ab\'ert that he recently proved
  that every infinite profinite group is nice \cite{Ab}.} 

\bt\lab{t:main2}
Suppose that every infinite profinite group is nice.\footnotemark[\value{footnote}] Then every
nondiscrete locally compact group in which every open 
subgroup is of index at most $\cof(\iN)$ is also nice; that is, there
exists a compact set $C$ of Haar measure zero such that the group can be
covered by $\cof(\iN)$ many left translates of $C$.
\et

Similarly to Corollary \ref{c:char} we also have the following.

\bcor
Suppose that every infinite profinite group is nice.
Then it is consistent with the axioms of set theory that a locally compact
group $G$ can 
be covered by less than $2^\om$ many translates of a compact nullset iff $G$
is nondiscrete and has no open subgroup of index at least $2^\om$. 
\ecor

The main goal of this section is to prove Theorem \ref{t:main2}. We start with
the Lie case. We use \cite{MZ} as the main reference, so note that
Lie groups are \emph{not} assumed to be second countable.

\bt\lab{t:Lie}
Every nondiscrete Lie group in which the identity component has index
at most $\cof(\iN)$ is nice; that is, it can be covered by
$\cof(\iN)$ many left translates of a compact set of Haar measure zero.
\et

\bp
Let $G$ be a Lie group as in the theorem. We can clearly assume that
$G$ is connected. 
Every compact neighbourhood of $e$ (the identity of $G$) generates an
open $\sigma$-compact subgroup, moreover, every open subgroup is actually
clopen. As $G$ is connected, 
we obtain that $G$ is
$\sigma$-compact, hence it possesses the Lindel\"of property. 
Therefore it suffices to show that there is a
neighbourhood of the identity that can be covered by
$\cof(\iN)$ many left translates of a compact set of Haar measure
zero.

Every nondiscrete Lie group contains one-parameter subgroups; that is; a
continuous homomorphic (not necessarily closed) image of $\RR$, see e.g.~\cite
[2.22.]{MZ}. Let $H\su G$ be the closure of such a subgroup, then 
$H$ is a closed connected commutative subgroup of $G$. 
By \cite[4.11.]{MZ} each closed subgroup of a Lie group is itself a Lie group, 
and so $H$ is actually a submanifold. If $G=H$,
we can apply Theorem \ref{main}, so we can assume that $H$ is a proper
subgroup. Let $M$ be a submanifold transversal to $H$ so that
$dim(H)+dim(M)=dim(G)$ and
\begin{equation}
\lab{HM}
H\cap M = \{e\}.
\end{equation}

\bl
There is a compact set $K \subset M$ which is a neighbourhood of $e$ (in M), so
that if $m: H \times K \rightarrow G$ is the restriction of the
multiplication map then
\begin{enumerate}[(i)] 
\item\lab{1} $m(H \times K) = HK$ is a neighbourhood of $e$
\item\lab{2} $m: H \times K \rightarrow HK$ is a homeomorphism.
\end{enumerate}
\el

\bp
It is well-known that if we use
the exponential map as a chart then the derivative of the $G\times G \to G$
multiplication map attains the form $(x,y) \mapsto x+y$ in the tangent spaces \cite{Wa}.
This implies that the
derivative of $m$ is nonsingular at $(e,e)$. Hence by the inverse function
theorem $m$ is a diffeomorphism in a neighbourhood of $(e,e)$. More precisely,
there exist open neighbourhoods $U, V$ and $W$ of $e$ in $H, M$ and $G$,
respectively, so that the restriction of $m$ is a smooth bijection of $U
\times V$ onto $UV = W$. 

This shows that (\ref{1}) holds for any choice of $K$ that is a
neighbourhood of $e$. 

Now we claim that
\begin{equation}
\lab{HW}
H \cap W \subset U.
\end{equation}
Indeed, if $h=uv \in H \cap UV$ then $v = u^{-1} h$. As $U\subset H$ we
obtain $v\in H$, and as $V\subset M$ by (\ref{HM}) we get $v=e$, so $h=u$.

Choose a compact neighbourhood $K\subset V$ of $e$ in $M$ so that
\begin{equation}\lab{K}
K K^{-1} \subset W.
\end{equation}

We claim that $K$ satisfies (\ref{2}), which will finish the proof of the lemma.

First we show that the map $m : H\times K\to HK$ is injective.  If $h_1
k_1=h_2 k_2$ then $h_2^{-1} h_1 = k_2 k_1^{-1}$.  Set $h=h_2^{-1} h_1 = k_2
k_1^{-1}$, then clearly  $h\in H$, and by (\ref{K}) we also have $h\in W$, hence
by (\ref{HW}) we obtain $h\in U$. 

Now we apply the fact that $m$ is a bijection between $U\times K$ and $UK$ (as
$K\su V$) to
the elements $hk_1=ek_2$. Indeed, $h, e \in U$ and $k_1, k_2\in K$, so
$k_1=k_2$ and $h_1=h_2$, proving that $M$ is injective.


Finally, we show that the inverse of $m$ is also
continuous. We use again the fact that $m$ is a smooth bijection between
$U\times K$ and $UK$.
So let $U_1, K_1$ be neighbourhoods of some $h\in H$ and $k\in
K$, respectively, then $h^{-1}U_1 \times K_1k^{-1}$ is a neighbourhood of
$(e,e)$. Hence its image $h^{-1}U_1 K_1k^{-1}$ contains a neighbourhood $W_1$
of $e$. Thus $hW_1k$ is a neighbourhood of $hk$ contained in the image of
$U_1\times K_1$ under $m$, proving that the inverse of $m$ is also
continuous.
\ep

Now we complete to proof of Theorem \ref{t:Lie}.

Fix a Haar measure $\mu_H$ on $H$, and consider a compact nullset $C$ in
$H$ as in Theorem \ref{main}.
The set $CK$ is compact, and $\cof({\iN})$ many left
translates of $CK$ cover $HK$, which is a neighbourhood of $e$ in
$G$. Therefore the proof of the theorem will be complete once we show the
following.

\bl\lab{null}
$CK$ is of $\mu_G$-measure zero, where $\mu_G$ is a left Haar measure on $G$.
\el

\bp
By the above lemma the multiplication map $H\times K \to HK$ is a
homeomorphism, hence $BK$ is Borel for every Borel set $B\su H$.
So we can define the set-function
\[
\mu: B \mapsto \mu_G(BK), \ \ (B\su H \textrm{ Borel}).
\]
It is easy to see that this is a left-invariant measure which is finite
for compact sets. We check that if $A\su H$ is a nonempty open (in $H$) set
then $\mu(A)>0$. Let $a\in A$, then $a^{-1} A$ is a
neighbourhood of $e$ in $H$. Clearly $\mu(A) = \mu_G(AK) = \mu_G(a^{-1}AK)>0$,
where the last expression is positive, since
$a^{-1}AK$ is a neighbourhood of $e$ in $G$.

By the uniqueness of Haar-measure \cite[17.B]{Ke}, there exists $c>0$, so that
$\mu=c\mu_H$, and so $\mu_G(CK)=\mu(C)=c\mu_H(C)=0$. This concludes the proof
of the lemma, and hence of the theorem.
\ep

\ep

\br
Lemma \ref{null} also follows from the construction of Haar measure via an invariant
smooth volume form, but we decided to use this alternative approach, which
establishes the lemma in a more direct fashion.

The proof of the above theorem with minor modifications shows that if $H$ is
a closed subgroup of a separable Lie group $G$, and $H$ can be covered by
$\kappa$ many left translates of a compact nullset then $G$ can also be
covered by $\kappa$ many left translates of a compact 
nullset. It would be interesting to see if this remains true in general, and
if so, if it could be used to establish our main theorem for profinite
groups.
\er

Next we consider the compact case. The following fact is most probably well
known. It was communicated to us by Ken Kunen.

\bs\lab{s:Kunen}
Every infinite compact group has a factor which is either an infinite
Lie group or an infinite profinite group.
\es

\bp
More precisely we show that if $G$ is an infinite compact group then
either it has an infinite Lie group factor or $G$ itself is profinite.

Denote by $U(n)$ the unitary group on $\CC^n$. 
By the Peter-Weyl theorem \cite[27.40]{HR} the set of all
representations of $G$ in the $U(n)$'s separate points of $G$,
hence $G$ is (isomorphic to) the inverse limit of the images of these
representations. If all these images are finite then $G$ is
profinite, and we are done. Otherwise $G$ has a factor that is an
infinite compact subgroup 
of some $U(n)$. But by \cite[4.11]{MZ} each closed subgroup of a Lie
group is itself a Lie group, so we are done.
\ep

Now we are ready to prove Theorem \ref{t:main2}. 

\bd
We say
that a topological group \emph{does not contain arbitrarily small
subgroups} 
if there is  a neighbourhood of the identity that contains no subgroup
apart from the one formed by the identity element itself. 

The identity
component of $G$ is denoted by $G_0$.
\ed

\bp (Theorem \ref{t:main2}) First note that if $H$ is a closed normal
subgroup in a topological group $G$ and every open subgroup of $G$ is
of index at most $\cof(\iN)$ then the same is true for $G/H$ and also
for every open subgroup of $G$.

Suppose that $G$ is a nondiscrete locally
compact group in which every open subgroup is of index at most
$\cof(\iN)$. We have to cover $G$ by at most
$\cof(\iN)$ many left translates of a compact left nullset.
By \cite[4.5. Cor.]{MZ} (actually, in every locally compact group $G$) 
there exists an open subgroup $G'
\su G$ and a compact normal subgroup $H$ of $G'$ such that $G'/H$ does
not contain arbitrarily small subgroups. 

$G'$ is clearly nondiscrete, since $G$ is nondiscrete and $G'$ is
open. As the index of $G'$ is at
most $\cof(\iN)$, it is sufficient to cover $G'$ by at most
$\cof(\iN)$ many left translates of a compact left nullset, hence we
can assume $G=G'$. 

So $H$ is a compact subgroup of $G$ such that $G/H$ does not contain
arbitrarily small subgroups. 

Now we separate two cases. First assume that $H$ is open. It suffices
to show that $H$ is nice. As above,
$H$ cannot be discrete, so it is infinite. By Statement \ref{s:Kunen}
either $H$ has 
an infinite profinite factor, in which case we are done by assumption
(and by Lemma \ref{l:factor}),
or $H$ has a factor which is an infinite Lie group. But an infinite
compact Lie group is clearly nondiscrete, and every open subgroup has
finite index, so we are done with this case by Theorem \ref{t:Lie}
(and again by Lemma \ref{l:factor}).

So we can assume that $H$ is not open, hence $G/H$ is not discrete.
By Lemma \ref{l:factor} it is sufficient to show that $G/H$ is nice.
By \cite[4.2.~Cor.~2]{MZ} if a locally
compact group does not contain arbitrarily small subgroups then the
identity component is open, hence $(G/H)_0$ is open in $G/H$. By
the remark at the beginning of the proof the index condition holds 
for $G/H$ too, so it is sufficient to show that $(G/H)_0$ is nice. 

As $G/H$ does not contain arbitrarily small
subgroups, the same holds for the subgroup $(G/H)_0$. 
By \cite[4.4.~Thm.]{MZ} a connected
locally compact group that does not contain arbitrarily small
subgroups is a Lie group, and clearly all these
requirements hold for $(G/H)_0$. Moreover, as $G/H$ is nondiscrete,
the same holds for the open subgroup $(G/H)_0$. 
Hence Theorem \ref{t:Lie} shows that $(G/H)_0$ is nice, finishing the
proof.
\ep

We conclude with some natural questions. 
Theorem \ref{t:main2} shows that the first two are equivalent.

\bq
Can we drop the assumption in Theorem \ref{main} that the group is
abelian? \!\!\footnote{We have been informed by M.~Ab\'ert that he recently
proved that every infinite profinite group is nice, which answers
these questions affirmatively \cite{Ab}.}
\eq

Or equivalently,

\bq
Suppose $G$ is an infinite profinite group. Is $G$ nice? That is, can
$G$ be covered by $\cof(\iN)$ many left translates of a compact set of
Haar measure zero? \!\!\footnotemark[\value{footnote}]
\eq

Of course in both questions it is also natural to replace $\cof(\iN)$ by
$<2^\omega$. In that case one can show that these
questions are also equivalent to the original Question 3.2 in \cite{ES}.

Our last question is a reformulation of \cite[Question 3.4]{ES}.

\bq
Suppose that $\kappa$ is a cardinal and $G_1, G_2$ are uncountable locally
compact (abelian) separable (Polish) groups such that $G_1$ can be covered 
by $\ka$ many translates of a suitably chosen compact nullset. Is the same 
true for $G_2$? 
\eq

\textbf{Acknowledgement.} The authors are indebted to 
Kenneth Kunen and Slawomir Solecki for some useful conversations. We are also
grateful to the anonymous referee for his useful suggestions.

\bigskip

\noindent
\textsc{M\'arton Elekes} 

\noindent
\textsc{Alfr\'ed R\'enyi Institute of Mathematics} 

\noindent
\textsc{Hungarian Academy of Sciences}

\noindent
\textsc{P.O. Box 127, H-1364 Budapest, Hungary}

\noindent
\textit{Email address}: \verb+emarci@renyi.hu+

\noindent
\textit{URL:} \verb+http://www.renyi.hu/~emarci+

\bigskip

\noindent
\textsc{\'Arp\'ad T\'oth} 

\noindent
\textsc{E\"otv\"os Lor\'and University, Department of Analysis}
 
\noindent
\textsc{P\'az\-m\'any P\'e\-ter s\'et\'any 1/c, H-1117, Budapest, Hungary}

\noindent
\textit{Email address}: \verb+toth@cs.elte.hu+

\bigskip

\end{document}